\documentclass[a4paper,twoside,10pt]{article}

\usepackage{amssymb}
\usepackage{amsmath}
\usepackage{epsfig}

\newtheorem{theorem}{Theorem}[section]
\newtheorem{lemma}[theorem]{Lemma}
\newtheorem{e-proposition}[theorem]{Proposition}
\newtheorem{corollary}[theorem]{Corollary}
\newtheorem{e-definition}[theorem]{Definition\rm}
\newtheorem{remark}{\it Remark\/}

\def\C{\mathbf C}
\def\D{\Delta}
\def\d{\delta}

\def\f{\phi}
\def\g{\gamma}

\def\l{\lambda}

\def\O{\mathbf O}
\def\p{\pi}

\def\R{\mathbf R}
\def\r{\rho}

\def\s{\sigma}
\def\t{\tau}

\def\Z{\mathbf Z}
\def\z{\zeta}

\title{Analytical invariants of quasi-ordinary hypersurface 
singularities associated to divisorial valuations}

\author{Pedro Daniel Gonz\'alez P\'erez, G\'erard Gonzalez-Sprinberg}

\date{}

\begin{document}

\maketitle

\begin{quote}\vskip 0.5\baselineskip
\noindent{\bf Abstract.} 
We study an analytically irreducible 
algebroid germ  $(X, 0)$ of complex singularity 
by considering  the filtrations of its analytic algebra, 
and their associated graded rings,
induced by the
{\it divisorial valuations}
associated to the irreducible components  
of the exceptional divisor of the normalized blow-up of
the normalization  $(\bar{X}, 0)$ of $(X, 0)$, centered at the point $0 \in \bar{X}$.
If  $(X, 0)$ is a quasi-ordinary hypersurface singularity, we obtain that 
the associated graded ring is a $\C$-algebra of finite type, 
namely the coordinate ring of a non necessarily normal affine toric variety
of the form $Z^\Gamma = \mbox{\rm Spec} \C [\Gamma]$, and we show that 
the semigroup  $\Gamma$ is an analytical  invariant of $(X, 0)$.
This  provides another 
proof of the analytical invariance of the {\it normalized 
characteristic monomials} of $(X, 0)$ (see \cite{Gau}, \cite{Lipman2} 
and  \cite{PP3}).
If  $(X, 0)$ is the algebroid germ of non necessarily normal toric variety, 
we  apply the same method to prove a local version of the isomorphism 
problem for algebroid germs of non necessarily normal toric varieties
(see \cite{Gu}, for the algebraic case).

\end{quote}

\section*{Introduction}

We study a class of algebroid hypersurface singularities
called {\it quasi-ordinary}. 
These singularities arise classically in Jung's approach 
to analyze surface singularities by using a finite projection to a smooth surface (see \cite{Jung}, \cite{Walker}).
Quasi-ordinary hypersurface singularities
are parametrized by {\it  quasi-ordinary branches}, 
certain class of fractional power series 
in several variables having a finite set of {\it distinguished} or {\it characteristic monomials},
 generalizing the {\it characteristic pairs} associated to a plane branch 
(see \cite{Mo}). 
These monomials determine many features of the geometry and topology of the singularity, 
for instance in the analytically irreducible case 
they define a complete invariant of the embedded topological type 
(see Lipman's and Gau's works \cite{Lipman2} and \cite{Gau}).
This characterization implies the analytical invariance of the  characteristic
monomials suitably {\it normalized} by an inversion lemma of Lipman
(see \cite{Gau}, Appendix and \cite{GP1}).
In the case of algebroid surfaces 
the analytical invariance of the normalized characteristic monomials was proved 
by Lipman and also by  Luengo, by building canonical sequences of monoidal and quadratic 
transformations which desingularize the germ, however these methods do not admit generalizations 
to quasi-ordinary hypersurfaces of dimension $\geq 3$
(see \cite{Lipman1} and \cite{Luengo}). 

In what follows $(S, 0)$ denotes the algebroid germ of an analytically irreducible 
quasi-ordinary hypersurface. 
In \cite{Tesis}, the first author introduced a semigroup $\Gamma$ associated 
to a fixed quasi-ordinary branch parametrizing the germ $(S, 0)$
and showed, using the characterization of the topological type, that 
different quasi-ordinary branches parametrizing $(S, 0)$
provide isomorphic semigroups. 
In this Note we prove that the 
semigroup $\Gamma$ is an analytical invariant
of the algebroid germ $(S, 0)$, without passing by Lipman's and Gau's
topological approach.
The isomorphism class of the  semigroup $\Gamma$  determines the normalized 
characteristic monomials  of the germ $(S, 0)$ 
(see \cite{Tesis}, \cite{GP1} or \cite{PP2}). 
Our result provides a proof of the analytical invariance 
of the normalized characteristic monomials of the
quasi-ordinary hypersurface $(S, 0)$.

Popescu-Pampu has given 
another proof of the analytical invariance of the 
  normalized characteristic monomials of the quasi-ordinary hypersurface $(S, 0)$ 
(see  \cite{PP3} and  \cite{PP2}). 
In his work \cite{PP3}, 
 which 
has a  flavor similar to that of the surface case \cite{PP2},
the structure of the singular locus of  $(S, 0)$ 
is used to build a sequence of blow-ups of $(\C^d, 0)$ (where $d = \dim S$), in terms  of
the presentation of the normalization of $(S, 0)$ as a quotient singularity.
This sequence is used to define a semigroup $\Gamma'$ depending only on the germ $(S, 0)$.
He proves,  using a fixed quasi-ordinary branch parametrizing  $(S, 0)$, that 
the semigroup  $\Gamma'$ is a linear projection of the semigroup $\Gamma$, 
eventually, of rank less than 
the rank of $\Gamma$. 
The normalized characteristic monomials are recovered 
from the semigroup $\Gamma'$ and from analysis of the equisingularity class of 
certain plane sections 
of $(S, 0)$, this latter technique is also  essential in Gau's work 
 \cite{Gau}.

Instead of defining an intrinsic semigroup associated to the singularity,
we study certain filtrations of its  analytic algebra, 
and their associated graded rings. This approach is inspired by the description, due to Lejeune-Jalabert, 
of the semigroup of a plane branch in terms of the graded ring associated 
to the filtration of its local ring  by the powers of the 
maximal ideal of its integral closure (see \cite{Mono} page 161).
If $(X, 0)$ is an analytically irreducible algebroid germ so is its normalization $(\bar{X}, 0)$.
We can study the singularity $(X, 0)$ by considering  the filtrations, and their associated graded rings,
induced by the
{\it divisorial valuations}
associated to the irreducible components  
 of the exceptional divisor of the normalized blow-up of $(\bar{X}, 0)$ centered at the point $0 \in \bar{X}$.
In general graded rings associated to divisorial valuations on normal or regular local rings are not 
necessarily Noetherian (see examples by  Cutkosky and Srinivas \cite{CS}, page 557,  and by
Cossart, Galindo and Piltant in \cite{CGP}).
We apply this strategy successfully  in two particular instances:

\begin{itemize}

\item
If  $(X, 0)$ is the algebroid germ of a non necessarily normal affine toric variety of the form, 
$Z^{\Lambda}= \mbox{\rm Spec } \C [ \Lambda ]$
at its zero dimensional orbit $0$,  we obtain that the  graded ring associated 
to $(X, 0)$ by any of these divisorial valuations, is  
a $\C$-algebra of finite type, namely the coordinate ring of the affine toric variety $Z^{\Lambda}$.
We apply a theorem of Gubeladze on the {\it isomorphism problem} for commutative monoid rings 
(see \cite{Gu}) to recover 
the semigroup $\Lambda$ from the toric variety $Z^\Lambda$.
This means that the semigroup $\Lambda $ is an analytic invariant of the algebroid germ $ (X, 0)$,
and it does not depend on the choice of  
toric structure on $(X,0)$. 
A topological proof of this result, in the case of  normal simplicial toric singularities, 
has been given by  Popescu-Pampu in \cite{PP4}.

\item
If $(X, 0)$ is the  germ of quasi-ordinary hypersurface $(S, 0)$ 
we obtain that the graded ring associated to $(S, 0)$ by any of these divisorial valuations, is  
a $\C$-algebra of finite type, namely the coordinate ring of the affine toric variety 
$Z^\Gamma: = \mbox{\rm Spec} \, \C [ \Gamma ]$.  The algebroid germ $(Z^\Gamma, 0)$  of the affine toric variety 
$Z^\Gamma$
at its zero 
dimensional orbit $0$,  is of  geometric significance to the hypersurface germ $(S, 0)$, 
for instance both germs have the same normalization (see \cite{GP3}), and 
any of these divisorial valuations associates to $(S, 0)$ and  $ (Z^\Gamma, 0) $ 
a pair of isomorphic graded rings. 
The analytical invariance of the semigroup $\Gamma$ follows simultaneously for the algebroid germs 
 $(S, 0)$ and  $ (Z^\Gamma, 0) $.
The same strategy provides a proof of the analytical 
invariance of the semigroup associated to 
a {\it toric quasi-ordinary hypersurface}.
This kind of singularity
belongs to  certain class of  
branched coverings of
a germ of normal affine toric variety at the special point,
which are unramified over the torus, and 
plays an important role in the toric embedded resolutions of 
quasi-ordinary hypersurface singularities described in \cite{GP3}. 
It is shown in \cite{GP3} that the germs $(S, 0)$ and $(Z^\Gamma, 0)$ 
have simultaneous toric embedded resolutions, which are described from the properties
of the semigroup $\Gamma$ by a method inspired by 
that of Goldin and Teissier for plane branches (see  \cite{Rebeca}).
\end{itemize}
 
\section{A reminder of toric geometry}

We give some definitions and notations (see \cite{Oda}, \cite{Oda2}, \cite{Ew}, \cite{Fulton}
and \cite{TE} for proofs).
If $N$ is a lattice we denote by $M$ the dual lattice, by
$N_\R$ the real vector space spanned by $N$.
A {\it rational convex polyhedral cone} $\sigma$  in  $N_\R$, a {\it cone} in what follows, 
is the set of  non
negative linear combinations of vectors $a^1, \dots, a^s  \in N$.
The cone $\sigma$ is {\it strictly convex} if
$\sigma$ contains no linear subspace of dimension $>0$.
We denote by $\stackrel{\circ}{\sigma}$ the { \it relative interior} of a 
cone $\sigma$.
The {\it dual} cone  $\sigma^\vee$ (resp. {\it orthogonal} cone 
$\sigma^\bot$) of $\sigma$ is the set
$ \{ w  \in M_\R / \langle w, u \rangle \geq 0,$  (resp. $ \langle w, 
u \rangle = 0$)  $ \; \forall u \in \sigma \}$.
A {\it fan\index{fan}} $\Sigma$ is a family of strictly convex
  cones  in $N_\R$
such that any face of such a cone is in the family and
the intersection of any two of them is a face of each.
The {\it support\index{support of a fan}} of the fan $\Sigma$ is the set
$ \bigcup_{\sigma \in \Sigma} \sigma \subset  N_\R$.

  A non necessarily normal affine toric
  variety is of the form
$Z^{\Lambda} = \makebox{Spec} \, \C [ \Lambda ]$ where $ \Lambda  $ is a
sub-semigroup of finite type of a lattice $ \Lambda + (- \Lambda)$ which it generates as a group.
The torus $Z^{ \Lambda + (- \Lambda)}$ is an open dense subset of $Z^\Lambda$, 
which acts on  $Z^\Lambda$, and this action 
extends
the action of the torus on itself by multiplication.
The semigroup $\Lambda$ spans the cone $\R_{\geq 0} \Lambda$ thus
we have  an inclusion of semigroups 
$\Lambda \rightarrow \bar{\Lambda} :=  \R_{\geq 0} \Lambda  \cap  ( \Lambda + (- \Lambda) ) $, defining an 
associated toric modification $Z^{\bar{\Lambda}} \rightarrow Z^{{\Lambda}}$, which  is the {\it normalization map}.
The cone
$\R_{\geq 0} \Lambda$ 
has a vertex if and only if 
there exists a zero dimensional orbit, and in this case this orbit is reduced to the point of $ 0 \in Z^{\Lambda}$ 
defined by the maximal ideal
$ \mathfrak{m}_\Lambda := (\Lambda - \{ 0 \} ) \C [ \Lambda ]$. 
The ring $\C [[ \Lambda ]] $ is the completion of the local 
ring of germs of holomorphic functions at $(Z^\Lambda, 0)$ with respect to 
its maximal ideal. 
 (See \cite{GKZ} for non necessarily normal toric varieties).
In particular, if  $\sigma$ is a cone in the fan $\Sigma$, the semigroup 
$\sigma^\vee  \cap M$  is of finite type,
it spans the lattice $M$
and the variety $Z^{\sigma^\vee \cap M}$, which
we denote also by  $Z_{\sigma, N}$ or by $Z_\sigma$ when the lattice is clear from the context, is normal.
If $\sigma \subset \sigma'$ are cones in the fan $\Sigma$ then we have an open immersion $Z_\sigma \subset Z_{\sigma'} $;
the affine varieties $Z_\sigma$  corresponding to cones in a fan $\Sigma$  
glue up to define the {\it toric variety\index{toric variety}}
$ Z_\Sigma$. The torus   $(  \C^*)^{ \mbox{rk}\, N} $ is the open dense subset $Z_{ \{ 0 \} } \subset  Z_\sigma$, and 
acts on $ Z_\sigma$ for each $\s \in \Sigma $; these actions
paste into an  action on $ Z_\Sigma$.
We say that a fan $\Sigma'$ is a {\it subdivision\index{fan subdivision}} of the fan
$\Sigma$ 
if both fans have the same support and if every cone of
$\Sigma'$ is contained in a cone of $\Sigma$.
The  subdivision $\Sigma'$ defines the {\it toric
modification\index{modification}}
$ \pi_{\Sigma'} : Z_{\Sigma '} \rightarrow   Z_\Sigma$
which is equivariant and 
induces an isomorphism between the tori.

We introduce, for each $\sigma \in \Sigma$, the   closed subset 
$\O_{\sigma}$  of $Z_\sigma$ defined
by the ideal $(X^w/ w  \in (\sigma^\vee - \sigma^\bot) \cap M)$ 
of $ \C [\sigma^\vee  \cap M ]$.
The coordinate ring of $\O_\sigma$ is $ \C [ \sigma^\bot \cap M]$.
The map that applies a cone $\sigma$ in the fan  $\Sigma$  to the set
$\O_\sigma \subset Z_\Sigma$ is a bijection between the relative interiors of the cones of the fan and
the orbits of the torus action which inverses inclusions of the closures.
The set  $\O_\sigma $ is the orbit of the {\it 
special point}   $o_\sigma $
defined by $X^u (o_\sigma) = 1 $ for all $u \in \sigma^\bot \cap  M $.
We have that $\dim \O_\s = \mbox{codim}\, \s$, in particular if $\t$ is an edge of $\Sigma$ 
 the closure $D(\t)$ of the orbit $\O_\t$ in $Z_\Sigma$ 
is a divisor.

If $\dim \s = \mbox{rk}\, N$, the orbit $\O_{\sigma}$ is reduced to the 
special point $o_\s$ of the affine toric variety $Z_\s$. 
The toric variety defined by the fan formed by the faces of the cone $\s$ coincides with 
the affine toric variety $Z_\s$. If $\Sigma$ is a subdivision of $\s$
we have that the exceptional fiber of the modification $ \pi_{\Sigma} : Z_{\Sigma } \rightarrow   Z_\s$ is equal to 
$\p_\Sigma ^{-1} (o_\s) = \bigcup_{\t \in \Sigma,\stackrel{\circ}{\t} \subset 
\stackrel{\circ}{\sigma}} \O_\t$ (see Proposition page 199, \cite{GS-LJ}).
Any non empty set $I \subset \s^\vee \cap M$ defines an integral polyhedron in $M_\R$
as the convex hull of the set $\bigcup_{a \in I} a + \s^\vee$.
We denote this polyhedron by 
${\mathcal N}_\s (I)$ or by ${\mathcal N} (I)$
if the cone $\s $ is clear from the context.
The face  of ${\mathcal N} (I)$ determined by $\eta \in \s$ 
is  the set ${\mathcal F}_\eta := \{  v \in {\mathcal N} (I) /  \langle \eta , v \rangle   = \inf_{v' \in {\mathcal N} (I)} \langle \eta , v' \rangle \}$.
All faces of ${\mathcal N} (I)$ are of this form, the compact faces are 
defined by vectors in $\stackrel{\circ}{\sigma}$.
The {\it dual Newton diagram} $\Sigma (I)$ associated to ${\mathcal N} (I)$
is the subdivision of $\s$ formed by the cones 
$\s( {\mathcal F} ) := \{  \eta \in \s \; /
\langle \eta , v \rangle   = \inf_{v' \in {\mathcal N} (I)} \langle \eta , v'
\rangle ,  \; \forall v \in {\mathcal F}
\}$, for ${\mathcal F}$ running through the faces of ${\mathcal N} (I)$.
If $\Sigma = \Sigma (I) $,  the modification 
$\p_\Sigma : Z_\Sigma \rightarrow Z_\s$ is the {\it normalized blowing up} of 
$Z_\s$ centered at the monomial ideal defined by $I$ in  $\C [ \s^\vee \cap M ]$
(see \cite{TE}, Chapter I, section 2). 
The support of a series $\f = \sum c_u X^u $ in $\C [[ \s^\vee \cap M ]] $ is the set of 
exponents of monomials with non zero coefficient. 
The {\it Newton polyhedron} ${\mathcal N}(\f)$ is the integral  polyhedron 
associated to the support of $\f$. 
If $\eta \in \s$  we denote by $\f_{| \eta} $ the {\it symbolic restriction} 
$\sum_{u \in {\mathcal F} \cap M } c_u X^u$of $\f$ to the face
${\mathcal F}$ determined by $\eta$.

\section{Quasi-ordinary hypersurface singularities} 

A germ of algebroid hypersurface  $(S,0 ) \subset  (\C^{d+1}, 0) $ 
is 
{\it quasi-ordinary} 
if there exists a finite morphism $(S, 0) \rightarrow (\C^d,0)$
(called a {\it quasi-ordinary projection\index{quasi-ordinary projection}}) such that the discriminant locus is 
contained (germ-wise) in a normal crossing divisor.
In suitable coordinates depending on this projection, 
the hypersurface $(S, 0)$ has an equation $f=0$, 
where $f \in \C [[ X ]] [Y]$ is 
a {\it quasi-ordinary polynomial}: a  Weierstrass polynomial 
with discriminant $\D_Y f$ of the form
$\D_Y f = X^\d \epsilon$, where $\epsilon$ is a unit in the ring 
$ \C [[ X ]]$ of formal power series in the variables 
$X= (X_1, \dots, X_d)$ and $\d \in \Z^d_{\geq 0}$.
We will suppose from now on that the germ $(S, 0)$ is analytically irreducible at the origin, i.e., the 
polynomial $f$ is irreducible. 
The Jung-Abhyankar theorem 
guarantees that the roots of $f$, the associated quasi-ordinary branches,  
are fractional power series in the ring $\C [[ X^{1/m} ]]$ for some $m \in \Z_{\geq 0}$,
(see \cite{Abhyankar}). 
If the series 
$\{ \z ^{(l)} \}_{l =1}^{\deg f}  \subset \C [[ X^{1/m} ]]$ are 
the roots of $f$, the discriminant of $f$ is equal to
$\D_Y f = \prod_{i\ne j} (\z^{(i)} - \z^{(j)})$
hence 
each factor $\z^{(t)} - \z^{(r)}$
is of the form a monomial times a unit in $\C [[ X^{1/m} ]]$.
These monomials (resp. their exponents) are called {\it characteristic} or
{\it distinguished}.
The characteristic exponents $\l_1, \dots, \l_g $ of $f$ 
can be labeled 
in such a way that  $\l_1\leq \dots \leq \l_g$ coordinate-wise (see \cite{Lipman2}).
Since the polynomial $f$ is irreducible, 
we can  identify the analytic algebra of the germ $(S, 0)$ 
with the ring  
$ \C [[  X ]] [ \z] $, for $ \z$ any fixed  quasi-ordinary branch $ \z$ parametrizing $(S, 0)$.

We denote the lattice $\Z^{d}$ by $M$ and  by $M_j$ the lattice
$\Z^{d} + \sum_{\l_i < \l_{j+1} } \Z \l_i  $,
for $j=1, \dots, g$ with the convention $\l_{g+1} = + \infty$;
the index $n_j$  of the lattice $M_{j-1} $ in $M_j$ non zero,  for $j=1, \dots,g$
(see \cite{Tesis}, \cite{GP3} and \cite{Lipman2}).
We denote by $N$ (resp. by $N_j$) the dual lattice of $M$ (resp. of $M_j$, for $j=1, \dots,g$)
and by   $\r \subset N_\R$ the cone spanned by the
dual basis of  the canonical basis of $M$.
With these notations we have  $\C [[ X ]] = \C [[ \r^\vee \cap M ]]$.
Since the exponents of the quasi-ordinary branch $\z$ belong to the semigroup
$\r^\vee \cap M_g$ we obtain a 
ring extension 
$\C [[  \r^\vee \cap M ]] [\z] \rightarrow \C [[ \r^\vee \cap M_g ]]$
which is 
the inclusion in the integral closure  (see Proposition 14, \cite{GP3}).
Geometrically, the germ of toric variety $ (Z_{\r, N_g}, o_\r)$ is
the normalization $(\bar{S}, 0)$  of $(S, 0)$; more generally 
the normalization of a germ of quasi-ordinary singularity, non necessarily hypersurface, 
is a toric singularity (see \cite{PP4}).
In our case,  the monomial ideal $\mathfrak{m}_{\r^\vee \cap M_g}: =(\r^\vee \cap M_g - \{ 0 \})
\C [[\r^\vee \cap M_g]] $ is the maximal ideal 
of the closed point $0 \in \bar{S}$.

\section{The invariance of the semigroups}

Let $(X,0)$ be an analytically irreducible algebroid germ. Denote by $p : (\bar{X}, 0)  \rightarrow (X, 0)$ the 
normalization map and  by $b: (B, E)  \rightarrow (\bar{X}, 0)$
the normalized blow up centered at  $0 \in \bar{X}$.
If $D$ is any irreducible component of  the exceptional divisor $E$ of 
the modification $b$ we denote by $\nu_D$
the associated divisorial valuation of the field of fractions 
$K$ of the analytic algebra $R$ of the germ $(X, 0)$. 
We have that $\nu_D (h)$ is equal to the vanishing order of $h\circ b$ 
at the component $D$,  for $h \in K - \{ 0\}$.
The valuation $\nu_D$ defines  a filtration of 
$R$ with ideals 
$\mathfrak{p}_k  = \{ \f \in R / \nu_D ( \f) \geq k \}$ for $k \geq 0$.
We denote by $\mbox{\rm gr}_D (X,0)$ the associated graded ring
$\mbox{\rm gr}_D (X,0): = \bigoplus_{k \in \Z_{\geq 0} } \mathfrak{p}_k / \mathfrak{p}_{k+1}$.
The graded ideal $\mathfrak{m}_D (X,0)  := \bigoplus_{k \in \Z_{\geq 1} } \mathfrak{p}_k / \mathfrak{p}_{k+1}$ is 
maximal and the pair $(\mbox{\rm gr}_D (X,0),\mathfrak{m}_D (X,0))$ 
depends only on the analytic algebra $R$ and the exceptional divisor $D$.

Suppose that $(\bar{X}, 0)$ is the germ of an affine toric variety $Z_{\r, N}$ at its zero dimensional orbit
(which is assumed to exists). 
The normalized blow up of the germ $(\bar{X}, 0) = (Z_{\r, N}, o_\r)$ centered at $0 = o_\r  $ is the toric modification
$\p_\Sigma : Z_\Sigma \rightarrow Z_{\r, N} $ 
where $\Sigma = \Sigma ( \r^\vee \cap M - \{ 0 \} )$. 
By definition, the modification $\p_\Sigma$ is an isomorphism outside the 
origin $o_\r \in Z_{\r, N}$ hence the  exceptional divisor $E$ is 
equal to the exceptional fiber $\p_\Sigma^{-1} (o_\r)$. 
The irreducible components of this divisor
are of the form $D(\t) $, where $\t$ runs through the edges of $\Sigma$
with $\stackrel{\circ}{\t} \subset \stackrel{\circ}{\r} $.
There exists at least one  such edge $\t$, we fix it and we denote $D(\t)$ by $D$
and by $n$ the integral lattice vector of $N$ in the edge  $\t$.  
The following property of the divisorial valuation associated to a toric divisor 
is used  
by Gonzalez-Sprinberg and Bouvier in the algebraic case (see \cite{Bouvier} and \cite{B-GS}). 
\begin{lemma}
If 
$0 \ne  \f = \sum c_u X^u \in \C [[ \r^\vee \cap M ]]$ 
then 
$\nu_D (\f) = \min_{c_u \ne 0} \langle n, u \rangle$.
\end{lemma}
{\em Proof}. Let $v $ denote any vertex of the compact face 
${\mathcal F} $ defined by $n$ on the polyhedron ${\mathcal N} (\f)$, 
we have that $\f $ factors by 
$X^v$ in a neighborhood of the divisor $D$ in the open chart 
$Z_\t$ of $Z_\Sigma$.
We have that 
$X^{-v} \f = X^{-v} \f_{| {\mathcal F}} + r $ and
the exponents of the terms in $r$  belong to $\t^\vee - \t^\bot$, 
thus $r$ vanish on $\O_\t$ and then on $D$.
We have also that $\nu_D  ( X^{-v} \f_{| {\mathcal F}}) =  0$
since   $ X^{-v} \f_{| {\mathcal F}}$ is a polynomial in $\C [\t^\vee \cap M]$
with non zero constant term.
We deduce from these facts that 
$\nu_D (\f) = \nu_D( \f_{| {\mathcal F}} )  = \nu_D( X^v )$
which is equal to $\langle n, v \rangle = \min_{c_u \ne 0} \langle n, u \rangle$
(see \cite{Fulton}, page 61). $\Box$

\subsection{The toric case}

Suppose that $(X,0)$ is the germ of a non necessarily normal affine toric variety $Z^{\Lambda}$ 
at its zero orbit $0$ (which is assumed to exist). The normalization $Z^{\bar{\Lambda}}$ 
is the affine toric variety $Z^{\bar{\Lambda}} = Z_{\r, N} $,  where $N$ is 
the dual lattice of $\Lambda + (- \Lambda)$ and $\r \subset N_\R$ is the dual cone of $ \R_{\geq 0} \Lambda$. 
We consider the 
graduation $\C [ \Lambda ] ^{(n ) } $ of the algebra $\C[   \Lambda ]$, with homogeneous terms 
$H_k  := \bigoplus_{u \in \Lambda, \langle n, u \rangle = k} \C X^u$, for 
$k \in \Z_{\geq 0}$ induced by the primitive vector $n \in \stackrel{\circ}{\t}$ associated to 
the exceptional divisor $D$.
We denote by $\mathfrak{m}_\Lambda ^{(n) } $ the maximal graded ideal $ \mathfrak{m}_\Lambda ^{(n) }:=
 \bigoplus_{k \geq 1} H_k  $ and by $\mathfrak{m}_\Lambda =  (\Lambda -\{ 0 \} ) \C [\Lambda]  $ the same ideal without the graded structure.

\begin{e-proposition} \label{t}
The pair 
$(\mbox{\rm gr}_D (Z^\Lambda,0),\mathfrak{m}_D (Z^\Lambda,0))$
is isomorphic to $(\C [\Lambda] ^{(n ) },\mathfrak{m}_\Lambda ^{(n ) } ) )$.
\end{e-proposition}
{\em Proof}. 
The analytic algebra   of the algebroid germ $(X,0)$ is isomorphic
to  $\C$-algebra $\C[[ \Lambda ]]$.
It is sufficient to prove that the map $\mathfrak{p}_k / \mathfrak{p}_{k+1} \rightarrow H_k$ 
given by $\f \in \C [[ \Lambda ]]$, $ \f \in \mathfrak{p}_{k} $,  
\begin{equation} \label{D}
 0 \ne \f \; \mbox{mod} \;  \mathfrak{p}_{k+1} \mapsto \f_{|n},  \mbox{ and } 
0 \; \mbox{mod} \;  \mathfrak{p}_{k+1} \mapsto 0,  
\end{equation}
is well defined and extends to a  graded isomorphism 
$ \mbox{\rm gr}_D (Z^\Lambda,0) \rightarrow \C [ \Lambda ]^{(n)}$.

We have that  if $0 \ne \f \in \C [[ \Lambda ]] $ and if $\eta \in \stackrel{\circ}{\r}$
the symbolic restriction $\f_{|\eta}$ belongs to $\C [ \Lambda]$, conversely  
given any $u \in \Lambda$ there exists 
$\f = X^u \in   \C [[  \Lambda ]]$ such that $\f_{|\eta} = X^u$.
It follows that 
$H_k - \{ 0 \}$ is equal to $\{ \f_{|n}  /  \f \in \mathfrak{p}_k | \mathfrak{p}_{k+1} \}$ (where 
$\mathfrak{p}_k | \mathfrak{p}_{k+1} = \{ \f \in \mathfrak{p}_k /  \f  \notin \mathfrak{p}_{k+1} \}$) .
If $0 \ne \f, \f' \in  \mathfrak{p}_{k} | \mathfrak{p}_{k+1}$ we have that
$ \f_{|n} = {\f'}_{|n}$ if and only if $\f = \f' 
 \; \mbox{mod} \;  \mathfrak{p}_{k+1}$
since the terms which may differ on $\f $ and $\f '$, 
viewed in 
$\C [[ \Lambda]]$,  have exponents of the form,  $u$ with
$\langle n, u \rangle > k$. We obtain that 
the map $\mathfrak{p}_k / \mathfrak{p}_{k+1} \rightarrow H_k$
is an isomorphism of vector spaces over $\C= H_0 = \mathfrak{p}_0 / \mathfrak{p}_{1}$.
The face of the Newton polyhedron ${\mathcal N} (\f \psi)$ defined by $n$ 
is equal to
the Minkowski sum of the faces defined by $n$ on the polyhedra 
 ${\mathcal N} (\f)  $ and ${\mathcal N} ( \psi) $ (see \cite{Ew}, page 105).
We deduce from this that 
$(\f \psi)_{| n} = \f_{| n}  \psi_{| n}$
and therefore that the map
$ \mbox{\rm gr}_D ( Z^\Lambda,0) \rightarrow \C [ \Lambda ]^{(n)}$
is a graded isomorphism.
$\Box$

We deduce from this result a local version of the {\it isomorphism problem} 
for germs of toric varieties (see \cite{Gu} for the affine case).
A topological proof of the following corollary
has been given by  Popescu-Pampu in \cite{PP4} 
in the case of  normal simplicial toric singularities.

\begin{corollary} \label{c}
Let $Z^{\Lambda}$ and $Z^{\Lambda'}$ be two non necessarily normal affine toric varieties
with zero orbits $0$ and $0'$ respectively.
If there exists an isomorphism of algebroid germs $(Z^{\Lambda},0) \cong (Z^{\Lambda'}, 0') $
the semigroups $\Lambda$ and $\Lambda'$ are isomorphic. 
\end{corollary}
{\em Proof}.
For any component $D$ of the exceptional divisor $E$ the pair 
obtained from  $(\mbox{\rm gr}_D (Z^\Lambda,0),\mathfrak{m}_D(Z^\Lambda,0))$
by forgetting the graduation is isomorphic to $(\C [\Lambda],\mathfrak{m}_\Lambda)$, thus it 
does not depend on the divisor $D$.
We obtain  isomorphic pairs of $\C$-algebra and maximal ideal
$(\C [ \Lambda ] ,\mathfrak{m}_\Lambda )$ and $(\C [ \Lambda' ] ,\mathfrak{m}_{\Lambda'} )$.
We apply  
Gubeladze's theorem 2.1 case (a) \cite{Gu},  to deduce the
existence of an isomorphism  of semigroups $\Lambda \cong \Lambda '$. $\Box$

\subsection{The quasi-ordinary hypersurface case}

We suppose now that  $(X, 0)$ is the  germ of quasi-ordinary hypersurface $(S,0)$.
In \cite{Tesis}
the semigroup 
$\Gamma := \Z^d_{\geq 0} + \g_1 \Z_{\geq 0} +  \cdots + \g_g \Z_{\geq 0} \subset \r^\vee \cap M_g $
where ${\g}_{1} =  \l_1$ and ${\g}_{j+1} = n_j {\g}_{j} +  \l_{j+1} -  \l_{j}$
for $j= 1, \dots, g-1$, is 
associated to a fixed quasi-ordinary branch $\z$ parametrizing 
the quasi-ordinary hypersurface germ $(S,0)$
(see \cite{GP1})
following the analogy to the case of plane curves (see \cite{Mo}).
The zero dimensional orbit $0$ of the affine toric variety $Z^\Gamma$
corresponds to the maximal ideal $\mathfrak{m}_\Gamma :=(\Gamma -\{ 0 \})\C [\Gamma]$.
In \cite{Tesis} it is shown that the normalizations of $(Z^\Gamma,0)$ and $(S,0)$ coincide.
We have that  $(\bar{X},0) = (Z^{\bar{\Gamma}},0) = (Z_{\r, N_g}, o_\r) $, with the notations of section 1
and 2.

\begin{e-proposition} \label{p}
For any irreducible component $D$ of the exceptional divisor 
of the normalized blow up of $(\bar{X},0)$ centered at $0$,
the pairs of graded ring and maximal graded ideal
$(\mbox{\rm gr}_D (S,0),\mathfrak{m}_D (S,0))$
and $(\mbox{\rm gr}_D (Z^\Gamma,0),\mathfrak{m}_D (Z^\Gamma,0))$
are isomorphic.
\end{e-proposition}
{\em Proof}. 
An irreducible component $D$ of the exceptional divisor of the normalized blow up of $(\bar{X},0)$
centered at $0$, corresponds to a primitive integral vector $n \in N_g$ in the interior of the cone $\r$.
The analytic algebra of $(S,0)$ is isomorphic to $\C [[ \r^\vee \cap M]] [\z]$, for $\z$ a fixed quasi-ordinary branch
parametrizing $(S,0)$. We consider the graded ring $\C [ \Gamma ] ^{(n)}$.
By proposition   \ref{t}
it is sufficient to prove that the map $\mathfrak{p}_k / \mathfrak{p}_{k+1} \rightarrow H_k$ 
given by formula (\ref{D}) for  $\f \in \C [[ \r^\vee \cap M]] [\z]$, $ \f \in \mathfrak{p}_{k}$,
is well defined and extends to a  graded isomorphism 
$ \mbox{\rm gr}_D (S,0) \rightarrow \C [ \Gamma ]^{(n)}$.

We have that  if $0 \ne \f \in \C [[ \r^\vee \cap M]][\z] $ and if $\eta \in \stackrel{\circ}{\r}$
the symbolic restriction $\f_{|\eta}$ belongs to $\C [ \Gamma]$, conversely  
given any $u \in \Gamma$ there exists 
$\f \in   \C [[ \r^\vee \cap M]][\z] $ such that $\f_{|\eta} = X^u$
(see Proposition 2.16 \cite{Tesis} or Proposition 3.1 and Theorem 3.6 of \cite{GP1}).
This property allows us to proceed with proof of the statement exactly in the same way as in 
proposition   \ref{t}.
$\Box$

\begin{corollary} \label{H}
If $\z$ and $\z'$ are two quasi-ordinary branches parametrizing the germ 
$(S, 0)$ then there exists an isomorphism of the corresponding 
semigroups  $\Gamma $ and $\Gamma'$.
\end{corollary}
{\em Proof}. 
For any irreducible component $D$ of the exceptional divisor of the normalized blow-up
of $(\bar{S},0)$, 
the pair of $\C$-algebra and maximal ideal  
obtained from $(\mbox{\rm gr}_D (S,0),\mathfrak{m}_D (S,0))$
by forgetting the graduation does not depend on the component $D$.
Indeed, it 
is isomorphic to $(\C [\Gamma],\mathfrak{m}_\Gamma)$
and to $(\C [\Gamma'],\mathfrak{m}_{\Gamma'}) $ by propositions \ref{t} and \ref{p}.
We apply 
Gubeladze's theorem 2.1 case (a) \cite{Gu},  to deduce the
existence of an isomorphism  of semigroups $\Gamma \cong \Gamma '$. $\Box$

\begin{remark}
Proposition \ref{p} and corollary \ref{H} extend to the case of 
{\it toric quasi-ordinary hypersurfaces} in \cite{GP3}. 
The arguments used in the proof of proposition \ref{p} above 
correspond to Proposition 31 Section 4.1 and Proposition 34 
Section 4.2 in \cite{GP3}.
\end{remark}

\section*{Acknowledgements}The first author is grateful to B. Teissier and P. Popescu-Pampu 
for their useful discussions and  suggestions, and to {\em Institut 
Fourier (Grenoble)} and  to 
 {\em Universidad de La Laguna} 
for their hospitality. 
The first  author was supported by a Marie Curie Fellowship of the European Community
programme ``Improving Human Research Potential and the Socio-economic Knowledge Base'' 
under contract number HPMF-CT-2000-00877.

\noindent
Both authors are grateful to the {\it Tokyo Institute of Technology}, who provide hospitality 
for the Conference ``Convex bodies and toric geometry'' in honor of  Professor Tadao Oda, in which 
the results of this Note were presented.

{\small 
{\sc Pedro Daniel Gonz\'alez P\'erez,} \\
Universit\'e de Paris 7, 
Institut de Math\'ematiques, 
Case 7012;
2, Place Jussieu,
75005 Paris, France.\\
gonzalez@math.jussieu.fr
}

{\small 
{\sc G\'erard Gonzalez-Sprinberg},\\
Institut Fourier, UMR 5582 du CNRS, 
Universit\'e de Grenoble I, 
 BP 74, 
38402 
Saint-Martin d'H\`eres, France\\
gonsprin@ujf-grenoble.fr
}


{\small

\begin{thebibliography}{00}



\bibitem{Abhyankar}{\sc Abhyankar, S.S.}, On the ramification of algebraic functions.
{\em Amer. J. Math.},  {\bf 77}. (1955), 575-592.

\bibitem{B-GS}{\sc Bouvier C., Gonzalez-Sprinberg G.,} Syst\`eme g\'en\'erateur minimal, diviseurs essentiels et G-d\'esingularisations de vari\'et\'es toriques, {\em T\^ohoku Math. J. }, Volume {\bf 47} (1995),  125-149.

\bibitem{Bouvier}{\sc Bouvier C.}, Diviseurs essentiels, composantes essentielles des vari\'et\'es toriques 
singul\`eres, {\em Duke Math. J.} Volume {\bf 91}, No 3, (1998), 609-620.


\bibitem{CGP} {\sc Cossart, V., Galindo,  C. and Piltant, O.} Un exemple effectif de gradu\'e non noetherien associ\'e \`a une valuation divisorielle, {\em Ann. Inst. Fourier (Grenoble)}, {\bf 50}, 1 (2000), 105-112.

\bibitem{CS} {\sc Cutkosky, D. and Srinivas V.} On a problem of Zariski on dimensions of linear systems, {\em Ann. Math.}, {\bf 137}, (1993), 173-189.



\bibitem{Ew}{\sc  Ewald,  G., } 
{\em Combinatorial Convexity and Algebraic Geometry,}
Sprin\-ger-Verlag, 1996.

 \bibitem{Fulton}{\sc  Fulton, W., } 
 {\em Introduction to Toric Varieties,}
 Annals of Math. Studies (131), Princenton University Press, 1993.


\bibitem{Gau}{\sc  Gau,  Y-N.,} 
{\em Embedded Topological classification of quasi-ordinary singularities,} 
Memoirs of the American Mathematical Society  388, 1988.

\bibitem{GKZ}{\sc Gel'fand, I.M., Kapranov, M.M. et Zelevinsky, A.V.}, 
{\em Discriminants, Resultants and Multi-Dimensional Determinants}, Birkh\"auser, Boston, 1994.

\bibitem{Rebeca}{\sc  Goldin, R., Teissier, B.,} 
  Resolving singularities of plane analytic branches with one toric morphism,
  {\em Resolution of Singularities}.
  Progress in Mathematics No. 181, Birkh\"auser-Verlag, 2000,
 315-340.

\bibitem{GP1}{\sc Gonz\'alez P\'erez P.D.}, The semigroup of a quasi-ordinary hypersurface,
 {\em J. Inst. Math. Jussieu}, (2003), {\bf 2} (3), 383-399.



\bibitem{Tesis}{\sc Gonz\'alez P\'erez P.D.}, Quasi-ordinary singularities via toric geometry, 
{\em Tesis Doctoral}, Universidad de La Laguna, (2000).

\bibitem{GP3}{\sc Gonz\'alez P\'erez P.D.}, 
Toric embedded resolutions of quasi-ordinary hypersurface singularities,
to appear in {\em Ann. Inst. Fourier (Grenoble)}.



\bibitem{GS-LJ}
 {\sc  Gonzalez-Sprinberg, G.,  Lejeune-Jalabert, M.},  Mod\`eles canoniques plong\'es. I.  
{\em Kodai Math. J.} 14 (1991), no. 2, 194-209. 




\bibitem{Gu} {\sc Gubeladze J.}, The isomorphism problem for commutative monoid rings, 
{\em Journal of Pure and Applied Algebra.},  {\bf 129}. (1998), 35-65.

\bibitem{Jung}{\sc  Jung,  H.W.E.,} 
 Darstellung der Funktionen eines algebraischen K\"orpers zweier
 unabha\"angigen Ver\"anderlichen $x$, $y$ in der Umgebung einer stelle
 $x=a$, $y=b$, {\em  J.Reine Angew. Math.},  {\bf 133} (1908),
 289-314.

 
\bibitem{TE}{\sc Kempf, G., Knudsen, F., Mumford D. , St Donat, B.},{\em Toroidal    Embeddings},
 Springer Lecture Notes in Mathematics No 339, Springer Verlag 1973.



\bibitem{Lipman1}{\sc  Lipman,  J.,} 
 Quasi-ordinary singularities of surfaces in $\C^3$, 
{\em Proceedings of Symposia in Pure Mathematics}, Volume {\bf 40} (1983), Part 2, 161-172.

\bibitem{Lipman2}{\sc  Lipman,  J.,} 
{\em Topological invariants of quasi-ordinary singularities,}
Memoirs of the American Mathematical Society  388, 1988.



\bibitem{Luengo}{\sc  Luengo,  I.,} On the structure of embedded
   algebroid surfaces, 
{\em Proceedings of Symposia in Pure Mathematics}, Volume {\bf 40}
  (1983), 185-193.

 \bibitem{Oda2}{\sc  Oda, T., } 
{\em Torus embeddings and applications. Based on joint work with Katsuya Miyake}. Tata Institute of Fundamental Research Lectures on Mathematics and Physics. Springer-Verlag,  1978.


\bibitem{Oda}{\sc  Oda, T., } 
 {\em Convex Bodies and Algebraic Geometry,}
 Annals of Math. Studies (131), Springer-Verlag, 1988.

\bibitem{PP2}{\sc Popescu-Pampu, P.}, 
On the invariance of the semigroup of a 
quasi-ordinary surface singularity, {\em C.R. Acad. Sci. Paris, Ser. I  },  {\bf 334}, (2002), 1101-1106 

\bibitem{PP3}{\sc Popescu-Pampu, P.}, 
On the analytical invariance of the semigroups of a quasi-ordinary hypersurface singularity, 
to appear in {\em Duke Math. J.}

\bibitem{PP4}{\sc Popescu-Pampu, P.}, 
On higher dimensional Hirzebruch-Jung singularities, 
preprint 2003, available at http://arxiv.org, 
math.CV/0306118.

\bibitem{Mono}{\sc  Teissier, B.,} The monomial curve and its  deformations, Appendix in \cite{Mo}.


\bibitem{Walker}{\sc  Walker,  R.J.,} 
 Reduction of the Singularities of an Algebraic Surface, 
 {\em Annals of Maths.},   {\bf 36}, 2 , (1935), 336-365.



\bibitem{Mo}{\sc  Zariski, O.,}  {\em Le probl\`eme des modules pour les branches planes}, Hermann, Paris, 1986.

\end{thebibliography}

}
\end{document}